\documentclass{gtart_a}
\pdfoutput=1
\usepackage{pinlabel}

\title[Tight contact structures with trivial contact 
invariants]{Infinitely many universally tight contact manifolds with trivial
Ozsv\'ath--Szab\'o contact invariants}

\author{Paolo Ghiggini}
\givenname{Paolo}
\surname{Ghiggini}
\address{CIRGET\\
Universit\'e du Qu\'ebec \`a Montr\'eal\\
Case Postale 8888, succursale Centre-Ville\\\newline
Montr\'eal (Qu\'ebec) H3C 3P8\\
Canada}
\email{ghiggini@math.uqam.ca}
\urladdr{}

\volumenumber{10}
\issuenumber{}
\publicationyear{2006}
\papernumber{10}
\lognumber{0676}
\startpage{335}
\endpage{357}

\doi{10.2140/gt.2006.10.335}
\MR{}
\Zbl{}

\arxivreference{math.GT/0510574}
\arxivpassword{}   


\received{4 November 2005}
\revised{}
\accepted{26 December 2005}
\published{2 April 2006}
\publishedonline{2 April 2006}
\proposed{Peter Ozsv\'ath}
\seconded{Tomasz Mrowka, Ronald Stern}
\corresponding{}
\editor{}
\version{}

\AtBeginDocument{\let\bar\wbar}

\makeop{Spin}
\newcommand{\uline}[1]{{\mskip3mu\underline{\mskip-3mu #1 \mskip-5mu}\mskip5mu}}

\makeatletter
\def\cnewtheorem#1[#2]#3{\newtheorem{#1}{#3}[section]
\expandafter\let\csname c@#1\endcsname\c@thm}

\theoremstyle{plain}
\newtheorem{thm}{Theorem}[section]
\cnewtheorem{prop}[thm]{Proposition}
\cnewtheorem{lemma}[thm]{Lemma}
\cnewtheorem{cor}[thm]{Corollary}

\theoremstyle{definition}
\cnewtheorem{rem}[thm]{Remark} %\unnumbered{rem}
\cnewtheorem{dfn}[thm]{Definition}
\cnewtheorem{conj}[thm]{Conjecture}
\makeatletter

\begin{document}

\begin{asciiabstract}
In this article we present infinitely many 3-manifolds 
admitting infinitely many universally tight contact structures 
each with trivial Ozsvath-Szabo contact invariants. By known 
properties of these invariants the contact structures constructed 
here are non weakly symplectically fillable.
\end{asciiabstract}

\begin{abstract}
In this article we present infinitely many 3--manifolds 
admitting infinitely many universally tight contact structures 
each with trivial Ozsv\'ath--Szab\'o contact invariants. By known 
properties of these invariants the contact structures constructed 
here are non weakly symplectically fillable.
\end{abstract}

\maketitle

%%%%%%%%%%%%%%%%%%%%%%%%%%%%%%%%%%%%%%%%%%%%%%%%%%%%%%%%%%%%%%%%%%
\section{Introduction}
Recently  Ozsv{\'a}th and Szab{\'o} introduced a new isotopy invariant 
$c(\xi)$ for contact $3$--manifolds $(Y, \xi)$ belonging to the 
Heegaard Floer homology group $\widehat{HF}(- Y)$. They  
proved \cite{O-Sz:cont} that $c(\xi)=0$ if $\xi$ is an overtwisted contact structure, 
and that $c(\xi) \neq 0$  if $\xi$ is Stein fillable. 
Later, they introduced also a refined version of the contact 
invariant denoted by $\underline{c}(\xi)$ taking values in the 
so-called Heegaard Floer homology group 
with twisted coefficients. They proved \cite[Theorem~4.2]{O-Sz:genus} that 
$\underline{c}(\xi) \neq 0$ if $(Y, \xi)$ is weakly fillable.

The Ozsv\'ath--Szab\'o contact invariants have been successfully used 
to prove tightness for several manifolds which had resisted to any 
previously known technique: see 
the papers \cite{lisca-stipsicz:4,lisca-stipsicz:5,lisca-stipsicz:6}
by Lisca and Stipsicz. 
This fact raised the hope that these
invariants could be non trivial for any tight contact structure. 
However in \cite{ghiggini:3} we showed that the  untwisted 
contact invariant reduced modulo $2$ can vanish even for
weakly symplectically fillable contact structures. Those examples
however left open the question whether the twisted invariants were 
non trivial for every tight contact structures.  In this article 
we give the following negative answer.

\begin{figure}
\labellist\small
\pinlabel {$0$} [l] at 265 95
\pinlabel {$-\frac{1}{r_1}$} [t] at 35 43
\pinlabel {$-\frac{1}{r_2}$} [t] at 98 25
\pinlabel {$-\frac{1}{r_3}$} [t] at 163 25
\pinlabel {$-\frac{1}{r_4}$} [t] at 223 40
\endlabellist
\centerline{\includegraphics[width=5cm]{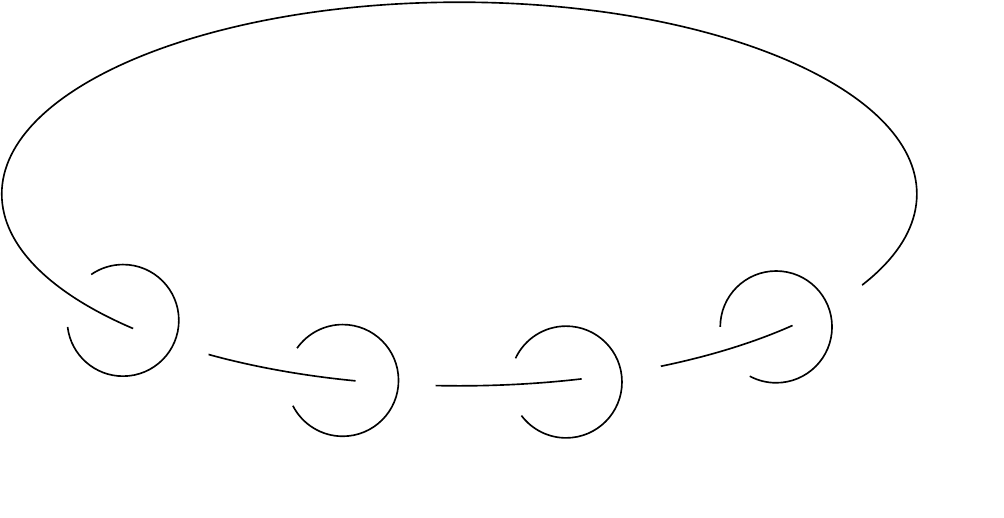}}
\label{definition.fig}
\caption{The surgery diagram of $M(r_1,r_2,r_3,r_4)$}
\end{figure}
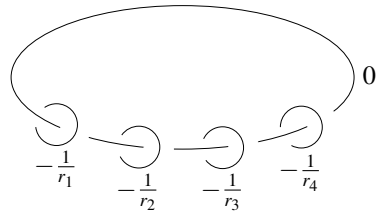

\begin{thm}\label{principale}
For any choice of coefficients $r_i \in (0,1) \cap \Q$ the Seifert 
manifold $M(r_1,r_2,r_3,r_4)$ defined by the surgery diagram in \fullref{definition.fig} admits infinitely many pairwise non 
isomorphic universally tight contact structures with trivial 
Ozsv\'ath--Szab\'o contact invariants.
\end{thm}

Although \fullref{principale} will be proved for the untwisted 
invariant, it holds for the twisted ones as well because 
$M(r_1,r_2,r_3,r_4)$ is a rational homology sphere for our choice of
Seifert coefficients, and for rational homology spheres 
the twisted and the untwisted contact invariants coincide. The 
following corollary is therefore a consequence of the 
non triviality of the twisted invariant for weakly symplectically 
fillable contact structures \cite[Theorem~4.2]{O-Sz:genus}, and of 
the fact that $H_1(M(r_1,r_2,r_3,r_4), \Q)=0$. 
\begin{cor}
For any choice of coefficients $r_i \in (0,1) \cap \Q$ the Seifert 
manifold $M(r_1,r_2,r_3,r_4)$ admits infinitely many pairwise non 
isomorphic universally tight contact structures which are not 
weakly symplectically fillable.
\end{cor}

This corollary provides the first universally tight contact 
structures which are known to be non weakly symplectically 
fillable, therefore it answers negatively to 
a question \cite[Question~4]{etnyre-ng:problems} of Etnyre and Ng, asking whether any 
universally tight contact $3$--manifold has to be weakly fillable.

\begin{rem}
The non fillable contact structures constructed here are all homotopic
to Stein fillable contact structures. Our construction contrasts with
all previous non fillability results, whose proofs relied on homotopic 
properties of the contact structures.
\end{rem} 

\subsection*{Acknowledgements}
This work was initiated while the author was visiting Princeton University
supported by the NSF Focused Research Grant FRG-024466. We would like to
thank Andr\'as Stipsicz for several useful discussions about Heegaard
Floer homology, and Steve Boyer for his help with the diffeomorphisms
of Seifert manifolds.  The author was partially supported by a CIRGET
fellowship and by the Chaire de Recherche du Canada en alg\`ebre,
combinatoire et informatique math\'ematique de l'UQAM.

%%%%%%%%%%%%%%%%%%%%%%%%%%%%%%%%%%%%%%%%%%%%%%%%%%%%%%%%%%%%%%%%%
\section{Construction of the contact structures}\label{costruzione}
In this section for any natural number $n$ we will construct a 
universally tight contact structure $\xi_n$ on $M(r_1, r_2, r_3, r_4)$. 

We denote by $M'(r_1, \ldots ,r_k)$ the Seifert manifold over $D^2$ with 
$k$ singular fibres with Seifert coefficients $r_1, \ldots ,r_k$.
We can decompose the manifold $M(r_1, r_2, r_3, r_4)$ as 
\[ M(r_1, r_2, r_3, r_4)= M'(r_1, r_2) \cup T^2 \times [-1,1] \cup M'(r_3, r_4). \]
The orientation on $T^2 \times \{ -1 \}$ is given by the inward normal 
convention, while the orientation on $T^2 \times \{ 1 \}$ is given by the 
outward normal convention, therefore $\partial M'(r_1, r_2)$ is identified 
to $T^2 \times \{ -1 \}$, and $T^2 \times \{ 1 \}$ is identified to $-\partial M'(r_3, r_4)$.

Since our construction will not use the Seifert coefficients of 
the fibres in any specific way, except for the fact $r_i \in (0,1)$, 
we will suppress them from the notation, and will call 
$M=M(r_1, r_2, r_3, r_4)$, $M_1'= M'(r_1, r_2)$, and $M_2'= M'(r_3, r_4)$.

By Hatcher \cite[Proposition~2.2 and Section~1.2]{hatcher}, $M_i'$ is a 
surface bundle over $S^1$. Let $\Sigma_i$ be the fibre, and $\phi_i$ be the 
monodromy of this bundle. The Seifert fibration on $M_i'$ 
restricted to $\Sigma_i$ is a 
branched cover $\Sigma_i \to D^2$ with finite fibre, and $\phi_i$ is a deck 
transformation, therefore $\phi_i$ has finite order $n_i$. 

Call ${\mathcal C}_i$ the set of all $1$--forms $\beta$ over $\Sigma_i$ such 
that
\begin{enumerate}
\item $d \beta$ is a volume form on $\Sigma_i$,
\item $\beta|_{\partial \Sigma_i}$ is a volume form on $\partial \Sigma_i$.
\end{enumerate}
By Thurston--Winkelnkemper \cite{thurston-winkelnkemper} the sets
${\mathcal C}_i$ are nonempty and convex. We define a $\phi_i$--invariant $1$--form 
$\wbar{\beta}_i$ on $\Sigma_i$ by averaging as follows: we pick any 
$1$--form $\beta_i \in {\mathcal C}_i$ and define
\[ \wbar{\beta}_i = \frac{1}{n_i} \sum_{k=0}^{n_i-1} (\phi_i^k)^* \beta. \] 
By the convexity of ${\mathcal C}_i$ we have $\wbar{\beta}_i  \in 
{\mathcal C}_i$.
If $t$ is the coordinate of $[0,1]$, for any $K > 0$ the $1$--form 
$dt+ K \wbar{\beta}_i$ is a contact form on $\Sigma_i \times [0,1]$ which gives 
a well defined contact form on $M_i'$. By Gray's Theorem the contact 
structures on $M_i'$ obtained from different choices of $\beta_i$ are isotopic,
while the actual value of $K$ has no relevance for our construction.
We denote by $\xi_+$ the kernel of $dt+ K \wbar{\beta}_i$, and by $\xi_-$ 
the kernel of $-(dt+ K \wbar{\beta}_i)$.
 The following lemma is a 
straightforward computation on the contact forms.
\begin{lemma}\label{conto}
$\partial M_i'$ is a prelagrangian torus with slope $s_i$ with respect to $\xi_+$ 
and $\xi_-$. We can chose $K$ so that $-2 < s_1 < -(r_1+r_2)$ and 
$-2 < s_2 < -(r_3+r_4)$. The Reeb vector fields of the contact forms 
$\pm (dt+ K \wbar{\beta}_i)$  are tangent to the fibres of the Seifert 
fibration of $M_i'$, and $\xi_+$ and $\xi_-$ are transverse to the Seifert 
fibration. 
\end{lemma}
In the following we will always assume that  $K$ has been chosen so 
that the inequalities in \fullref{conto} hold.

On $T^2 \times [-1, 1]$ we consider the contact structures
\[ \alpha_n(s_1, s_2) = \ker \big ( \cos (\varphi_n(t))dx+ \sin (\varphi_n(t))dy \big ) \]
for a smooth function $\varphi_n \colon [-1,1] \to \R$ such that 
\begin{enumerate}
\item $\varphi_n'(t)>0$ for any $t \in [-1,1]$,
\item $[\frac{\varphi_n(1)- \varphi_n(-1)}{\pi}]=n$,
\item $T^2 \times \{ -1 \}$ has slope $s_1$ and $T^2 \times \{ 1 \}$ has slope $-s_2$,
\item $\varphi_n(-1) \in (0, \frac{\pi}{2})$.
\end{enumerate} 
Condition (1) implies that $\alpha_n(s_1, s_2)$ is a contact structure, 
condition (2) implies that it has twisting $n \pi$ in the sense of
Honda \cite[Section 2.2.1]{honda:1}, and condition (4) is simply a 
normalisation condition.
 The set of functions satisfying these conditions is convex, 
therefore by Gray's Theorem the isotopy type of $\alpha_n(s_1, s_2)$ 
depends only on $n$, $s_1$ and $s_2$. 
 
\begin{dfn}
We define the contact structures $\xi_n$ on $M(r_1, r_2, r_3, r_4)$ so 
that they coincide with $\xi_+$ on $M_1'$, with $\xi_-$ on $M_2'$, and with 
$\alpha_{2n}(s_1, s_2)$ on $T^2 \times [-1,1]$, where $s_i$ is the boundary slope of 
$M_i'$.
\end{dfn}

Following Colin and Honda \cite{colin-honda:1} we say that a 
contact structure is {\em hypertight} if it can be defined by a 
contact form 
whose associated Reeb vector field has no contractible periodic 
orbits. By Hofer \cite[Theorem~1]{hofer:1} hypertight contact 
structures are tight.
\begin{thm} \label{hypertight}
The contact structures $\xi_n$ on $M(r_1, r_2, r_3, r_4)$ are hypertight 
for any $n \geq 0$.
\end{thm}
\begin{proof}
By \cite[Lemma~9.1]{colin-honda:1} we can isotope $\xi_+$ and 
$\xi_-$ relative to the boundary, so that they are defined by contact 
$1$--forms which glue to the contact form of $\alpha_{2n}(s_1,s_2)$ to give 
a globally defined contact form for $\xi_n$. Moreover, the isotopy 
can be chosen so that the Reeb vector fields of $\xi_+$ and $\xi_-$ remain
transverse to the fibrations over $S^1$ defined on $M_1'$ and $M_2'$. 

Any periodic orbit of the Reeb flow must be completely contained 
in one of the three pieces $M_1'$, $M_2'$ and $T^2 \times [-1,1]$ in which 
$M$ has been decomposed, because the Reeb vector field is tangent 
to $\partial M_1'$ and to $\partial M_2'$. Moreover, the inclusions $M_i' \hookrightarrow M$ and 
$T^2 \times [-1,1] \hookrightarrow M$ induce injective maps between the fundamental 
groups, therefore a periodic orbit of the Reeb flow is contractible
in $M$ if and only it is contractible in the piece it is contained 
in. This implies that there are no contractible periodic orbits of 
the Reeb flow in $M$, because the Reeb vector field is transverse 
to the $S^1$--fibrations in $M_i'$, and in $T^2 \times [-1,1]$ its closed
orbits are homotopically non trivial in the incompressible tori 
$T^2 \times \{ c \}$. 
\end{proof}
In order to prove universal tightness for $\xi_n$ we need the 
following lemma about the coverings of Seifert manifolds.
\begin{lemma}\label{rivestimenti}
Let $M$ be a Seifert manifold with base $B$, and denote its 
universal covering by $\wwtilde{M}$. If $K$ is a compact set of 
$\wwtilde{M}$, then there exists a finite covering $\wwbar{M}$
 of $M$ for which the projection $\wwtilde{M} \to \wwbar{M}$ 
is injective on $K$.
\end{lemma}
\begin{proof}
We have to consider only the case when the universal cover is not 
already a finite cover itself. Consider the orbifold structure 
induced on $B$ by the Seifert fibration on $M$. Either $M$ has 
at least three singular fibres, or the genus of $B$ (as a surface) 
is $g >0$ because we have assumed that the universal cover of $M$ 
is infinite. In these cases $B$ is a good orbifold, therefore 
there is a finite orbifold covering $B' \to B$ such that $B'$ has 
no singular points: see Scott \cite[Theorem~2.3 and Theorem~2.5]{scott}. 

We pull back the Seifert fibration of $M$ to $B'$ in order to 
obtain a Seifert manifold $M'$ which fibres over $B'$ and a finite 
covering $M' \to M$. The Seifert fibration on $M'$ has no singular 
fibres because $B'$ has no singular points, therefore $M'$ is a 
circle bundle over $B'$.  This concludes the proof because the 
lemma holds trivially for circle bundles over surfaces, and the 
composition of finite coverings is still a finite covering. 
\end{proof}

\begin{cor} \label{utight}
The contact structures $\xi_n$ on $M(r_1, r_2, r_3, r_4)$ are universally 
tight for all $n \geq 0$.
\end{cor}
\begin{proof}
Suppose by contradiction that the universal cover of $(M, \, \xi_n)$ 
contains an overtwisted disc. Since the overtwisted disc is 
compact, by \fullref{rivestimenti} $(M, \, \xi_n)$
has an overtwisted finite cover. This is a contradiction
because any finite cover of a hypertight contact manifold is 
hypertight again, and therefore it is tight. 
\end{proof}

%%%%%%%%%%%%%%%%%%%%%%%%%%%%%%%%%%%%%%%%%%%%%%%%%%%%%%%%%%%%%%%%%%%
\section{Decomposition of the tight contact structures} \label{tre}
Let $M$ be a Seifert manifold with base $B$ and $k$ singular fibres $F_1,
\ldots ,F_k$, and let $U_i$ be a standard neighbourhood of $F_i$, $i=1,
\ldots ,k$. Then $M \setminus \bigcup_{i=1}^k U_i$ can be identified
with $S \times S^1$ where $S$ is the surface obtained by removing $k$
disjoint discs centred at the images of the singular fibres from the
base $B$. This diffeomorphism determines identifications of $- \partial
(M \setminus U_i)$ with $\R^2 / \Z^2$ so that $\bigl( {1 \atop 0} \bigr)$
is the direction of the section $S \times
\{ 1 \}$ and $\bigl({0 \atop 1}\bigr)$
is the direction of the regular fibres. In order to fix one among the
infinitely many product structures on $M \setminus \bigcup_{i=1}^k
U_i$ we also require the meridian of each $U_i$ to have slope $-
\frac{\beta_i}{\alpha_i}$ in $- \partial (M \setminus V_i)$ with
$\frac{\beta_i}{\alpha_1} = r_1$.

We also choose an identification between $\partial U_i$ and $\R^2 /
\Z^2$ so that $\bigl({1 \atop 0}\bigr)$
is the direction of the meridian of $U_i$ and $\bigl({0 \atop 1}\bigr)$
is the direction of a longitude. Notice
that $\partial U_i$ and $- \partial (M \setminus U_i)$ coincide as sets,
but are identified with $\R^2 / \Z^2$ in different ways. We can choose the
longitude on $U_i$ so that these two identifications are related by gluing
matrices 
$A_i \co  \partial U_i \to - \partial (M \setminus U_i)$ 
given by
\[
A_i= \left ( \begin{matrix}
\alpha_i & \alpha_i' \\ 
- \beta_i & -\beta_i'  
\end{matrix} \right )
\]
with $\beta_i \alpha_i'-\alpha_i \beta_i'=1$ and $0< \alpha_i' < \alpha_i$.

\begin{lemma} \label{thickening}
Let $\xi$ be a tight contact structure on $M$, and assume that 
$(M, \, \xi)$ has  a Legendrian regular fibre $L$ with 
twisting number $0$ (possibly after isotopy), which means that its 
contact framing coincides with the framing determined by the 
fibration. Then for $i=1, \ldots ,k$ there exist tubular neighbourhoods
$V_i$ of the singular fibres $F_i$ so that $\partial (M \setminus V_i)$
is a convex torus in standard form with infinite slope.
\end{lemma}

\begin{proof}
Make the singular fibres $F_i$ Legendrian with very low twisting 
numbers $n_i<0$, and consider their standard neighbourhoods $U_i$ for 
$i=1, \ldots ,k$. Without loss of generality we can assume $L \cap
U_i=\emptyset$ for $i=1, \ldots ,k$. Let $A_i$ be a convex
vertical annulus between $L$ and a Legendrian ruling curve
of $\partial (M \setminus U_i)$. By the Imbalance Principle,
Honda~\cite[Proposition~3.17]{honda:1}, $A_i$ produces a bypass attached
to $- \partial (M \setminus U_i)$ along a vertical Legendrian ruling
curve.  Then using the bypass attachment Lemma \cite[Lemma~3.15]{honda:1}
we can thicken $U_i$ until we obtain a convex solid torus $V_i$ such
that $- \partial (M\setminus V_i)$ has infinite slope.
\end{proof}
We call $(M \setminus \bigcup V_i, \, \xi|_{M \setminus \bigcup V_i})$ the
{\em background} of $(M, \, \xi)$. Since the background of 
$(M, \, \xi)$ has infinite boundary slopes, by 
Honda \cite[Section~4.3]{honda:2} it is isomorphic 
to an $S^1$--invariant tight contact structure 
$(S \times S^1, \, \xi_{\Gamma_{S_0}})$. Here $\xi_{\Gamma_{S_0}}$ denotes the 
$S^1$--invariant contact structure on $S \times S^1$ inducing the 
dividing set $\Gamma_{S_0}$ on a convex $\# \Gamma$--minimising section 
$S_0 \subset S \times S^1$ with Legendrian boundary. By 
\cite[Proposition~4.4]{honda:2} the isotopy class of $\xi_{\Gamma_{S_0}}$ 
is completely determined by  $\Gamma_{S_0}$.

Let $c_1$ be the smallest number for which the torus $T^2 \times \{ c_1 \}$ 
in the contact manifold $(T^2 \times [-1,1], \, \alpha_{2n}(s_1,s_2))$ has 
infinite slope, and let $c_2$ be the biggest one. Also, let $c_1'$ 
be the first number for which $T^2 \times \{ c_1' \}$ has slope $-2$ and let
$c_2'$ be the last number for which $T^2 \times \{ c_2' \}$ has slope $2$. 
$c_1'$ and $c_2'$ exist because of \fullref{conto}.
Call 
\begin{align*}
& (M_{c_1'}, \, \wtilde{\xi}_+')= (M_1', \, \xi_+) \cup 
(T^2 \times [-1,c_1'], \, \alpha_{2n}(s_1,s_2)|_{T^2 \times [-1,c_1']}) \\
& (M_{c_2'}, \, \wtilde{\xi}_-')= (M_2', \, \xi_-) \cup 
(T^2 \times [c_2', 1], \, \alpha_{2n}(s_1,s_2)|_{T^2 \times [c_2',1]}) \\
& (M_{c_1}, \, \wtilde{\xi}_+) = (M_1', \, \xi_+) \cup 
(T^2 \times [-1,c_1], \, \alpha_{2n}(s_1,s_2)|_{T^2 \times [-1,c_1]}) \\
& (M_{c_2}, \, \wtilde{\xi}_-) = (M_2', \, \xi_-) \cup 
(T^2 \times [c_2, 1], \, \alpha_{2n}(s_1,s_2)|_{T^2 \times [c_2,1]}). 
\end{align*}

\begin{lemma} \label{novert}
$(M_{c_1'}, \, \wtilde{\xi}_+')$ and $(M_{c_2'}, \, \wtilde{\xi}_-')$
contain no Legendrian curves with twisting number $0$ isotopic to
regular fibres. 
\end{lemma}
\begin{proof}
We prove the lemma only for $(M_{c_1'}, \, \wtilde{\xi}_+')$ because 
the proof for $(M_{c_2'}, \, \wtilde{\xi}_-')$ is the same.

Fix a rational number $r_3' < -2$, and consider the matrix with 
integral entries
\[ A(r_3')=  \left ( \begin{matrix}
        \alpha  &       \alpha' \\
        - \beta     &   - \beta'
\end{matrix} \right ) \]
such that $r_3'= \frac{\beta}{\alpha}$, $\alpha' \beta - \alpha \beta'=1$, and $0 < \alpha' < \alpha$. 
Applying $A(r_3')^{-1}$ to $- \partial M_{c_1'}$ we obtain a prelagrangian 
torus with slope
\[- \frac{\beta + 2 \alpha}{\beta' + 2 \alpha'}= - \frac{\alpha}{\alpha' - \frac{1}{\beta+2 \alpha}} \]
which is greater than the slope of the Seifert fibration 
$- \frac{\alpha}{\alpha'}$ because $\frac{\beta}{\alpha}< -2$. 

Put polar coordinates $(\rho , \theta)$ on $\R^2$, and call
\[ D(\rho_0)= \{ (\rho, \theta) \in \R^2 : \rho \leq \rho_0 \}. \]
We can choose $\rho_0$ so that $\big ( D(\rho_0) \times S^1, \, \ker(dz - \rho^2 d \theta)
\big )$ has prelagrangian boundary and boundary slope 
$$- \frac{\beta + 2 \alpha}{\beta' + 2 \alpha'}.$$
If we glue the tight solid torus
$$\big ( D(\rho_0) \times S^1, \, \ker(dz - \rho^2 d \theta) \big )$$
to $- \partial M_{c_1'}$
by the map $A(r_3')$, we obtain a contact structure $\tau$ on 
$M(r_1, r_2, r_3')$. This contact structure is transverse to the 
Seifert fibration because the contact planes do not twist enough
to become tangent to the fibres. 

Transverse contact structures in Seifert manifolds are tight by
 Lisca--Mati\'c \cite[Corollary 2.2]{lisca-matic:transverse}, therefore 
 Wu \cite[Theorem 1.4]{wu:1} implies that $(M(r_1, r_2, r_3'), \, \tau)$ 
contains no Legendrian curve with twisting number $0$ isotopic to 
a regular fibre. Consequently, $(M_{c_1'}, \, \wtilde{\xi}_+')$ 
contains no such a curve either.
\end{proof}

\begin{lemma} \label{menouno}
For any $i=1,2,3,4$ there exists a tubular neighbourhood $V_i'$ of 
the singular fibre $F_i$ such that $V_1', V_2' \subset M_{c_1'}$, 
$V_3', V_4' \subset M_{c_2'}$, and $- \partial (M \setminus V_i')$ has slope $-1$. Moreover
there exist collars $C_1$ of $\partial M_{c_1}$ in $M_{c_1} \setminus (V_1' \cup V_2')$ and 
$C_2$ of $\partial M_{c_2}$ in $M_{c_2} \setminus (V_3' \cup V_4')$ such that $\partial (M_{c_1} \setminus C_1)$ 
and $\partial (M_{c_2} \setminus C_2)$ are convex tori with slope $-1$.
\end{lemma}
\begin{proof}
Again, we prove the lemma only for $M_{c_1'}$. We perturb 
$\partial M_{c_1'}$ so that it becomes a convex torus in 
standard form with vertical ruling, then we make the singular 
fibres $F_1$ and $F_2$ Legendrian with very low twisting numbers $k_1$
and $k_2$, and take standard neighbourhoods $U_i$ of $F_i$. The slopes 
of $- \partial ( M_{c_1'} \setminus (U_1 \cup U_2))$ are
$$
\begin{cases}
- \frac{k_1 \beta_1 + \beta_1'}{k_1 \alpha_1 + \alpha_1'} & \text{on the
component corresponding to } \partial U_1,\\
- \frac{k_2 \beta_2 + \beta_2'\vrule width0pt height10pt }{k_2 \alpha_2 + \alpha_2'}& \text{on
the component corresponding to } \partial U_2,\\
\qquad2 & \text{on } - \partial M_{c_1'}.
\end{cases}$$
We can make $\smash{- \frac{k_i \beta_i + \beta_i'}{k_i \alpha_i +
\alpha_i'}}$ arbitrarily close to $-r_i \in (-1,0)$
by making $k_i$ very large in absolute value. 

Take convex vertical annuli 
$A_i$ between a vertical Legendrian ruling curve of $\partial M_{c_1'}$ and 
a vertical Legendrian ruling curve of the component of 
$\partial ( M_{c_1'} \setminus (U_1 \cup U_2))$ corresponding to $\partial U_i$ for $i=1,2$. 
By the imbalance principle of Honda \cite[Proposition~3.17]{honda:1} $A_i$ 
carries a bypass on the side of $\partial U_i$. Attaching this bypass we can 
thicken $U_i$ and reduce the slope of the component of 
$- \partial ( M_{c_1'} \setminus (U_1 \cup U_2))$ corresponding to $\partial U_i$. We can repeat 
this procedure until we get neighbourhoods $V_1'$ and $V_2'$ of $F_1$
and $F_2$ such that the slopes of $- \partial (M_{c_1'} \setminus (V_1' \cup V_2'))$ are 
$-1$, $-1$, and $2$. Denote by $T_1$ and $T_2$ the components of 
$- \partial (M_{c_1'} \setminus (V_1' \cup V_2'))$ corresponding to $V_1'$ and $V_2'$ 
respectively. Let $B$ be a convex vertical annulus 
between vertical Legendrian ruling curves of $T_1$ and $T_2$. The 
dividing set of $B$ contains no boundary parallel dividing curves, 
otherwise there would be a bypass attached vertically to $T_1$ or 
$T_2$ by \cite[Proposition~3.18]{honda:1}, and we could use this 
bypass to produce a torus with infinite slope inside $M_{c_1'}$. 
This would be a contradiction
because there are no Legendrian curves with twisting number $0$ 
isotopic to regular fibres in $M_{c_1'}$ by \fullref{novert}. 
After cutting $M_{c_1'} \setminus (V_1' \cup V_2')$ along $B$ and rounding the 
edges, by the Edge rounding lemma \cite[Lemma~3.11]{honda:1} we 
get a convex torus $T_3$ with slope $-1$ isotopic to $\partial M_{c_1}$. 
The collar $C_1$ is bounded by $T_3$ and  $\partial M_{c_1}$.
\end{proof}

Now we determine the isotopy type of the contact structures 
$\xi_n|_{V_i}$. We recall that $D^2 \times S^1$ admits exactly two universally
tight contact structures with $\# \Gamma_{\partial D^2 \times S^1}=2$ if its boundary 
slope is lesser than $-1$, and exactly one its boundary slope is
$-1$; see Honda \cite[Proposition~5.1(2)]{honda:1}. 
It follows from the computation of the relative Euler class of
the universally tight contact structures on $T^2 \times [0,1]$
\cite[Proposition~5.1]{honda:1} and from the correspondence between
tight contact structures on $T^2 \times [0,1]$ and on $D^2 \times S^1$ 
\cite[Proposition~4.15]{honda:1} that  all
the basic slices in the basic slices decomposition of a universally
 tight contact structures have the same sign. We define the sign 
of a universally tight contact structure on $D^2 \times S^1$ as the sign 
of any basic slices in its decomposition.
\begin{prop} \label{tori}
The contact structures $\xi_n|_{V_i}$ are universally tight for 
$i=1,2,3,4$. They are positive for $i=1,2$ and negative for 
$i=3,4$.
\end{prop}
\begin{proof}
$(V_i, \, \xi_n|_{V_i})$ is universally tight because $(M, \, \xi_n)$ is 
universally tight and the inclusions $\iota_i \colon V_i \to M$ induce injective 
maps $(\iota_i)_* \colon \pi_1(V_i) \to \pi_1(M)$. We can assume without loss of 
generality that $V_i'$ is contained in $V_i$ for $i=1,2,3,4$, then
$V_i \setminus V_i'$ is the outermost basic slice in the decomposition of 
$V_i$. Since all basic slice in the decomposition of a universally
tight contact structure have the same sign, the sign of 
$\xi_n|_{V_i \setminus V_i'}$ determines the sign of $\xi_n|_{V_i}$. By 
Ghiggini--Lisca--Stipsicz \cite[Lemma~2.7]{gls:2} the signs of $V_1 \setminus V_1'$ and of $V_2 \setminus V_2'$ 
are the same as the sign of $C_1$, and the signs of $V_3 \setminus V_3'$ and 
of $V_4 \setminus V_4'$ are the same as the sign of $C_2$. In applying 
\cite[Lemma~2.7]{gls:2} we must notice that here the sign of $C_i$ 
is computed after orienting $\partial M_{c_1}$ by the outward normal 
convention, while in \cite[Lemma~2.7]{gls:2} all boundaries are 
oriented by the inward normal. By a direct check on $\alpha_{2n}(s_1,s_2)$ it 
easy to see that $C_1$ is positive and $C_2$ is negative.
\end{proof}  
\begin{lemma}  \label{background1}
For $i=1, \ldots ,4$, let $V_i$ the neighbourhood of the singular 
fibre  $F_1$ obtained by applying \fullref{thickening} to 
$(M_{c_1}, \, \wtilde{\xi}_+)$ and to 
$(M_{c_2}, \, \wtilde{\xi}_-)$. 
Then the dividing set on a convex $\# \Gamma$--minimising section of the 
backgrounds of $(M_{c_1}, \, \wtilde{\xi}_+)$ and 
$(M_{c_2}, \, \wtilde{\xi}_-)$ has no boundary 
parallel dividing curves.
\end{lemma}
\begin{proof}
We prove the lemma only for  
$(M_{c_1}, \, \wtilde{\xi}_+)$ because the proof for 
 $(M_{c_2}, \, \wtilde{\xi}_-)$ is the same. 

Let $C_0 = M_{c_1} \setminus (V_1' \cup V_2' \cup C_1)$.  Since 
$\xi_n|_{C_0}= \wtilde{\xi}_+|_{C_0}$ is a tight contact 
structure with boundary 
slopes $-1$ and without Legendrian curves with twisting number $0$ 
isotopic to fibres, 
and $(C_1, \, \wtilde{\xi}_+)$ is a basic slice with 
boundary slopes $-1$ and $\infty$, we apply \cite[Lemma~2.7]{gls:2} and 
conclude that the dividing set on a convex $\# \Gamma$--minimising 
section of $M_{c_1} \setminus (V_1 \cup V_2)$ has no boundary parallel dividing 
arcs.
\end{proof}

\begin{figure}\small \centering
\begin{tabular}{ccc}
\includegraphics[width=5cm]{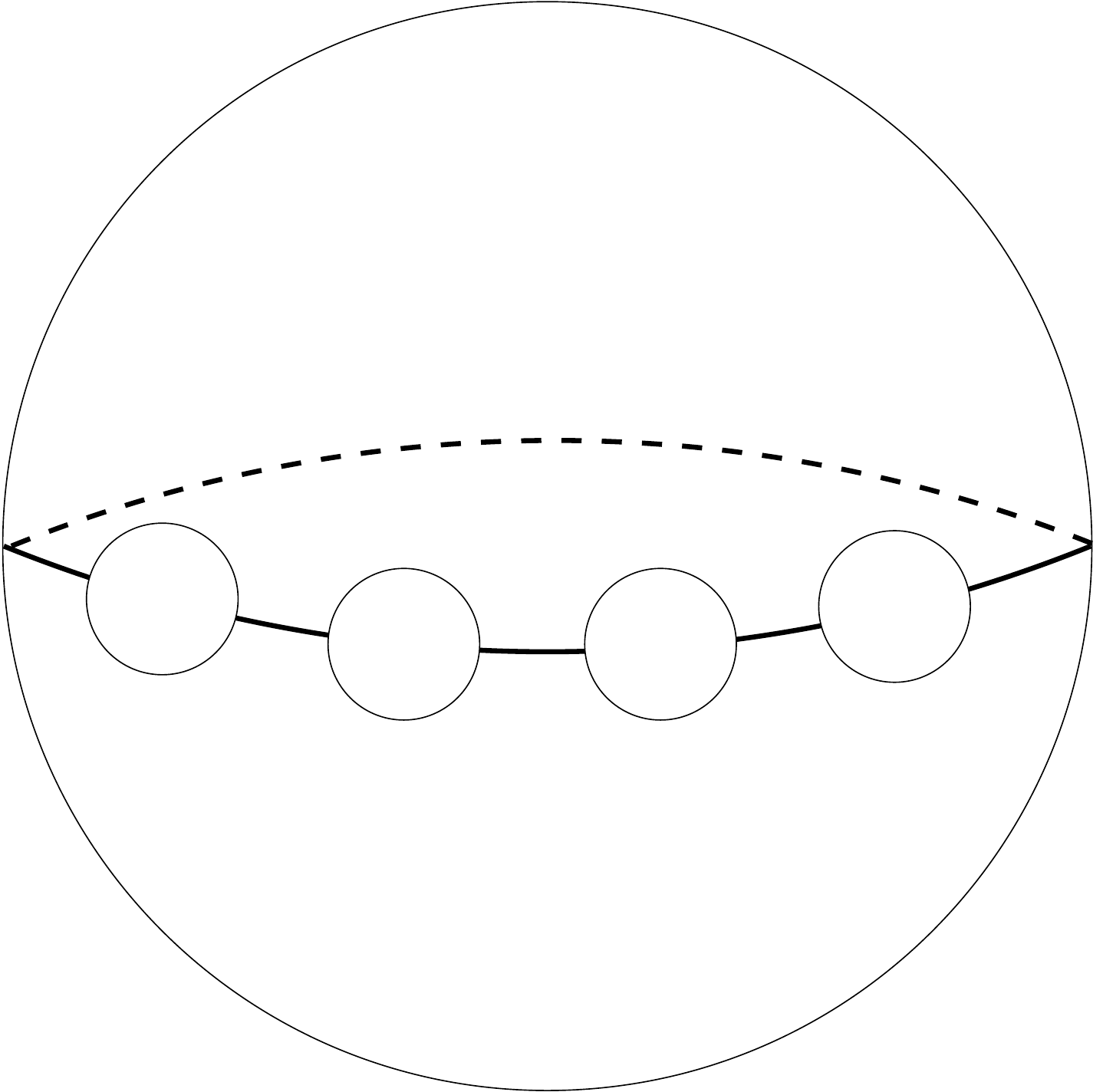}

& &

\includegraphics[width=5cm]{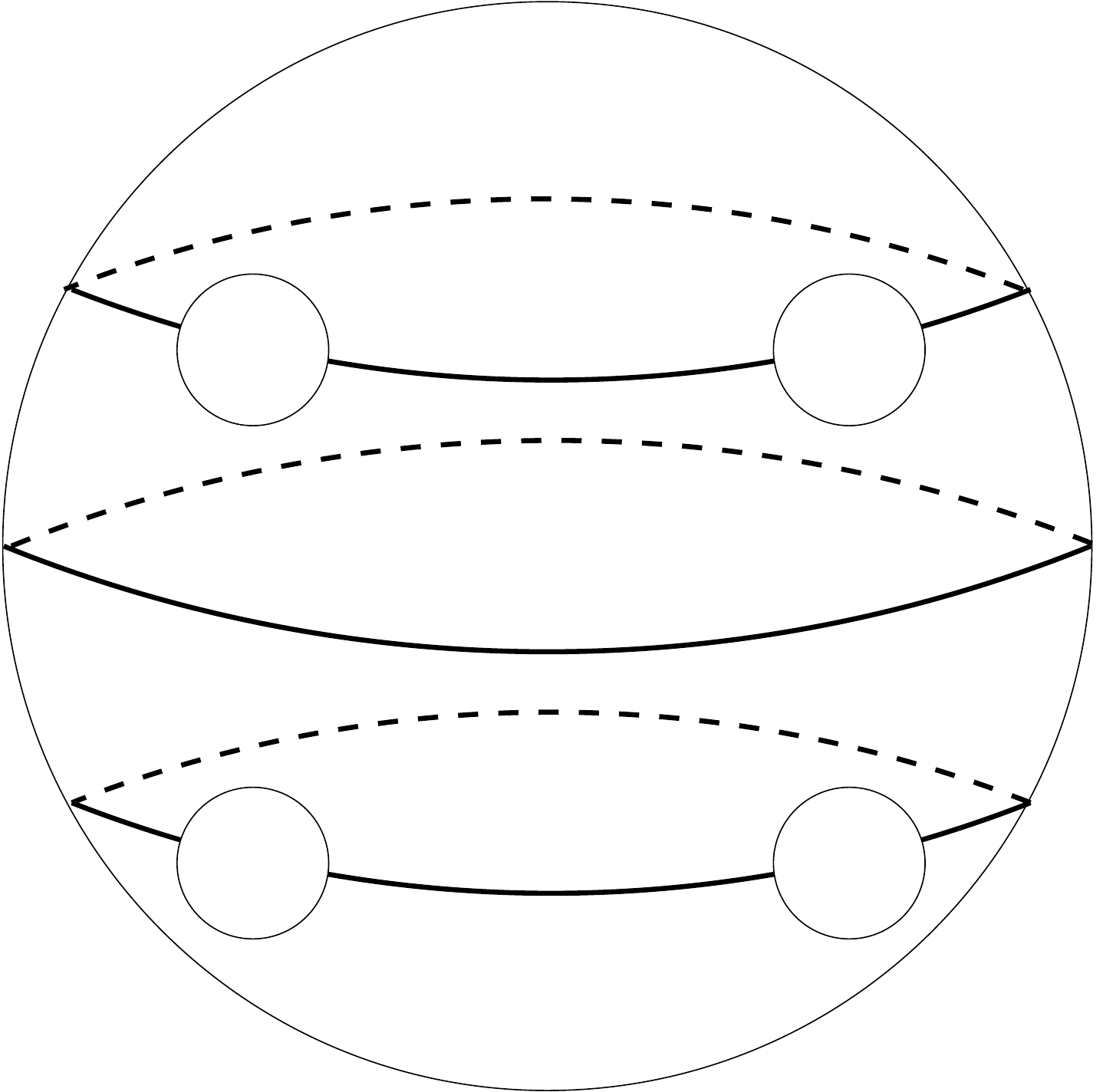} \\
(a) & & (b)

\end{tabular}
\caption{On the left the dividing set on a $\# \Gamma$--minimising 
section of the background of $(M, \, \xi_0)$. On the right the 
dividing set on a $\# \Gamma$--minimising section of the 
background of $(M, \,\xi_1)$.}
\label{background.fig}
\end{figure}

\begin{prop} \label{background2}
Let $S$ be a four-punctured sphere, and let $S_0 \subset S \times S^1$ be a 
convex $\# \Gamma$--minimising section in the background 
$(S \times S^1, \xi_{\Gamma_{S_0}})$ of $(M, \xi_n)$. Then, if we choose the 
neighbourhoods $V_i$ so that $V_1$, $V_2$ are contained in $M_{c_1}$ 
and $V_3$, $V_4$ are contained in $M_{c_2}$, the dividing set of $S_0$ 
is isotopic to one of the following:
\begin{enumerate}
\item four dividing arcs joining the components of $\partial S_0$ in 
cyclic order, as in a necklace, if $n=0$ (see \fullref{background.fig}(a)), or 
\item two dividing arcs joining $V_1$ to $V_2$, two dividing arcs 
joining $V_3$ to $V_4$, and $2n-1$ homotopically non trivial 
parallel dividing curves in-between (see \fullref{background.fig}(b)).
\end{enumerate}
\end{prop}  

\begin{figure}\centering
\includegraphics[width=5cm]{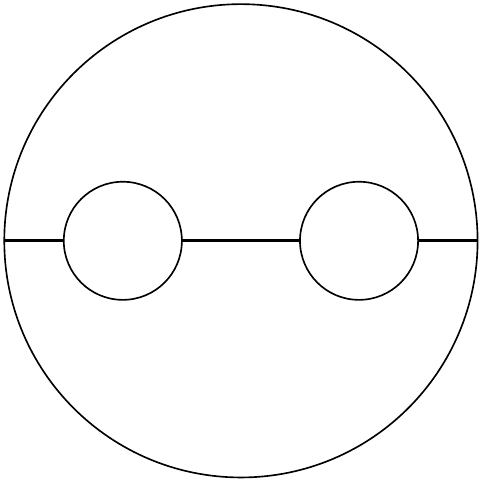}
\caption{The dividing set on $S_0^{(1)}$ and $S_0^{(3)}$}
\label{pantalone.fig}
\end{figure}

\begin{proof}
We divide $S_0$ into three pieces $S_0^{(1)} \subset M_{c_1}$, $S_0^{(2)} \subset T^2 \times [c_1, c_2]$,
and $S_0^{(3)} \subset M_{c_2}$, so that each piece is convex with Legendrian 
boundary and $\# \Gamma$--minimising in its relative homology class. Here 
we assume that the  product structure on $T^2 \times I$ has been deformed 
in small neighbourhoods of
$c_1$ and $c_2$ so that $T^2 \times \{ c_1 \}$ and $T^2 \times \{ c_2 \}$ has become 
convex tori with two dividing curves. If $n=0$ $c_1=c_2$, then we 
assume further that $c_1$ and $c_2$ have been replaced by $c_1- \epsilon$ and 
$c_2+ \epsilon$, so that $T^2 \times [c_1, c_2]$ has become the invariant 
neighbourhood of a convex torus. 

By \fullref{background1} the dividing set on $S_0^{(1)}$ and $S_0^{(3)}$ 
consists of three arcs connecting the boundary components in pairs,
as in \fullref{pantalone.fig}. 
If $n=0$ $T^2 \times [c_1,c_2]$ is an invariant neighbourhood of a convex 
torus, therefore the dividing set on  $S_0^{(2)}$ consists of two arcs 
connecting the two components of $\partial S_0^{(2)}$. 
If $n>0$, on the other hand, by 
Honda--Kazez--Mati\'c \cite[Proposition~2.2]{honda-kazez-matic:2} $S_0^{(2)}$ 
consists of a boundary parallel dividing arc for each
component of $\partial S_0^{(2)}$ and of $2n-1$ closed homotopically trivial 
curves. To see that the number of closed dividing curves is odd 
we glue the boundary components of $T^2 \times [c_1,c_2]$ together by the 
identity map, and observe from the equation defining $\alpha_{2n}(s_1,s_2)$
that we obtain a tight contact structure on $T^3$. This is possible
only if the boundary parallel arcs of $\smash{\Gamma_{S_0^{(2)}}}$
glue to give a
homotopically non trivial closed curve, and this happens only if the
number of closed dividing curves on $S_0^{(2)}$ is odd.
Gluing  $\smash{S_0^{(1)}}$, $\smash{S_0^{(2)}}$, and $\smash{S_0^{(3)}}$ together we obtain the 
dividing set on $S_0$. 
\end{proof} 

\section{Distinguishing the contact structures}
Our next goal is to prove that the contact manifolds 
$(M, \, \xi_n)$ are pairwise non isomorphic. We will follow the line
of Honda--Kazez--Mati\'c \cite[Section~4.3]{honda-kazez-matic:2}. For $i=1,2$ the fibre 
bundle $M_i' \to S^1$ with fibre $\Sigma_i$ defined in \fullref{costruzione}
  extends to a fibre bundle $M_{c_i} \to S^1$ with fibre 
$\wwbar{\Sigma}_i$ containing $\Sigma_i$. We define a surface with boundary 
$\wwhat{\Sigma} \subset M$ as follows. Identify $T^2 \times \{ c_i \}$ with $T^2$ for
$i=1,2$, and regard $\partial \wwbar{\Sigma}_i$ as a curve in $T^2$, then 
isotope $\wwbar{\Sigma}_1$ and $\wwbar{\Sigma}_2$ so that 
$\partial \wwbar{\Sigma}_1$ and $\partial \wwbar{\Sigma}_2$ minimise their geometric 
intersection and call $x_1, \ldots ,x_m$ the intersection points between 
$\partial \wwbar{\Sigma}_1$ and $\partial \wwbar{\Sigma}_2$ in $T^2$. Then for any 
intersection point $x_i$ join $\wwbar{\Sigma}_1$ and $\wwbar{\Sigma}_2$ by 
a small band thickening the segment $\{ x_i \} \times [c_1, c_2]$ so that the 
band intersects $T^2 \times \{ t \}$ in a linear arc whose slope is never 
vertical.  Call the resulting surface $\wwhat{\Sigma}$.

Now take two boundary-incompressible arcs $\gamma_1 \subset \wwbar{\Sigma}_1$ and 
$\gamma_2 \subset \wwbar{\Sigma}_2$ with endpoints on the same intersection 
points, and extend $\gamma_1 \cup \gamma_2$ over the bands to get a simple closed 
curve $\gamma \subset \wwhat{\Sigma}$. We define the framing of $\gamma$ to be the 
one coming from $\wwhat{\Sigma}$. Let ${\mathcal L}_n$ be the set of 
the Legendrian curves in $(M, \xi_n)$ which are smoothly isotopic to 
$\gamma$, and define the {\em maximal twisting} $t({\mathcal L}_n)$ to be 
the maximum attained by the twisting number of the curves in 
${\mathcal L}_n$.

\begin{prop} \label{torsione} $t({\mathcal L}_n) = -2n+1$. 
\end{prop}
\begin{proof}
This is a corollary of \cite[Proposition~4.9]{honda-kazez-matic:2}
once we have proved that the contact structures $\xi_n$ are isomorphic
to the contact structures $\zeta_{2k}$ defined in that article. To prove 
this we need to find a convex decomposition of $M_{c_1}$ and $M_{c_2}$
such that the dividing sets induced by $\wtilde{\xi}_+$ and 
$\wtilde{\xi}_-$ on the cutting surfaces are isotopic to the dividing 
sets induced by the contact structures constructed by Honda, Kazez and
Mati\'c \cite{honda-kazez-matic:2}. This means that we want 
the dividing set induced by $\wtilde{\xi}_+$ and $\wtilde{\xi}_-$ 
on any cutting surface to be boundary parallel. We will work out the 
details only for $(M_{c_1}, \, \wtilde{\xi}_+)$, because the
proof for  $(M_{c_2}, \, \wtilde{\xi}_-)$ is 
the same.

We take $\wwbar{\Sigma}_1$ as the first cutting surface; 
$M_{c_1} \setminus \wwbar{\Sigma}_1$ is diffeomorphic to $\wwbar{\Sigma}_1 \times 
[0,1]$, therefore we can further decompose it by cutting along 
discs of the form $\alpha_i \times I$, where 
$\{ \alpha_1, \ldots , \alpha_{1-\chi(\wwbar{\Sigma}_1)}\}$ is a set of 
properly embedded and pairwise disjoint arcs with boundary on 
$\partial \wwbar{\Sigma}_1$, such that $\wwbar{\Sigma}_1 \setminus (\alpha_1 \cup \ldots \cup 
\alpha_{1-\chi(\wwbar{\Sigma}_1)})$ is disc.
 
We can find a convex surface with Legendrian boundary isotopic
to $\wwbar{\Sigma}_1$ by ``patching  meridional discs'' of $V_1$ 
and $V_2$ as described in 
Etnyre--Honda \cite[proof of Proposition~3.5]{etnyre-honda:2} or
Ghiggini--Sch\"onenberger \cite[Section~4.1]{ghiggini-schonenberger}. Since the contact 
structures on  $V_1$ and $V_2$ are both universally tight and 
positive, the dividing sets on the meridional discs of $V_1$ and 
$V_2$ consist of boundary-parallel arcs cutting out  
regions all with the same signs, therefore when we patch the 
meridional discs their dividing arcs  join to 
give boundary-parallel dividing arcs in $\wwbar{\Sigma}_1$.

If we choose the 
arc $\alpha_i$ disjoint from the dividing set of 
$\wwbar{\Sigma}_i$ for all $i$, then $\partial (\alpha_i \times [0,1])$ intersects the dividing set
of $\partial (\wwbar{\Sigma}_1 \times [0,1])$ in exactly two points. When we
make $\alpha_i \times [0,1]$ convex with Legendrian boundary, its 
dividing set consists of exactly one arc, 
therefore we have obtained a convex decomposition of 
$(M_{c_1}, \, \wtilde{\xi}_+)$ of the required form.
\end{proof}

\begin{cor} \label{nonisot}
$\xi_n$ is not isotopic to $\xi_m$ if $n \neq m$
\end{cor}
\begin{proof}
$t({\mathcal L}_n)$ is clearly an isotopy invariant for $\xi_n$.
\end{proof}

\begin{thm} \label{nonisom}
$\xi_n$ is not isomorphic to $\xi_m$ if $n \neq m$
\end{thm}

\begin{proof}
Let $\phi \colon M \to M$ be a diffeomorphism of $M$ such that $\phi_*(\xi_n)=\xi_m$. 
Denote by $T_0$ a fibred torus in $M$ isotopic to $\partial M_1'$ and 
$\partial M_2'$. From the equation defining $\alpha_{2n}(s_1,s_2)$ and $\alpha_{2m}(s_1,s_2)$
it is immediate to check that for both contact structures the 
isotopy class of $T_0$ contains a prelagrangian torus. Any other 
incompressible torus in $M$ non isotopic to 
$T_0$ intersect $T_0$ persistently, therefore by 
Colin \cite[Theorem~1.6]{colin:1} any prelagrangian torus in $M$ is 
isotopic to $T_0$. This implies that $\phi(T_0)$ is isotopic to $T_0$. 
By Orlik \cite[Theorem~8.1.7]{orlik} $\phi$ is isotopic to a 
fibre-preserving diffeomorphism, therefore it defines 
an element $\wbar{\phi}$ in the mapping class group of the base 
orbifold $B$ of $M$. Call $C_0$ the projection of $T_0$ to $B$. Since
$\wbar{\phi}$ fixes $C_0$, it must be a product of Dehn twists 
around $C_0$. Let $C_1$ be an essential curve in $B$ which intersects
$C_0$ in exactly two points, and let $T_1$ be the pre-image of $C_1$ 
in $M$. $T_1$ is a fibred torus which splits $M$ in two submanifolds
$M_l$ and $M_r$, and we may also assume without loss of generality
that $F_1, F_3 \subset M_l$ and $F_2, F_4 \subset M_r$. Let $\phi'$ be a composition of 
$\phi$ with Dehn twists around $T_0$ and isotopies so that $\phi'(T_0)=T_0$ 
and $\phi'(T_1)=T_1$. Then we can assume that $T_0$ and $T_1$  are fixed 
not only as sets, but also pointwise  because the action of $\phi'$ 
on the homology of $T_0$ and $T_1$ must preserve the kernels of the 
maps 
\[ H_1(T_0) \to H_1(M_1'), \quad H_1(T_0) \to H_1(M_2') \]
and 
\[ H_1(T_1) \to H_1(M_l), \quad H_1(T_1) \to H_1(M_r). \]
Since these kernels are linearly independent, $\phi'$ acts trivial on
$H_1(T_0)$ and $H_1(T_1)$, therefore it is isotopic to the identity.

Since $M \setminus (T_0 \cup T_1)$ is a disjoint union of four solid tori and
$\phi'$ fixes their boundaries, $\phi'$ is isotopic relative to the 
boundary to the identity in any of the components of $M \setminus (T_0 \cup T_1)$,
 therefore it is isotopic to the identity on $M$. This implies that
$\phi$ is isotopic to a product of Dehn twists around $T_0$. We may 
further assume that
$\phi$ is supported in $M \setminus (M_1' \cup M_2') \cong T^2 \times [-1,1]$. The contact
structure $\alpha_{2n}(s_1,s_2)$ is invariant up to isotopy under Dehn 
twists around $T_0$, therefore $\xi_m= \phi_*(\xi_n)$ is isotopic to $\xi_n$. 
This proves that $m=n$, because $\xi_m$ is not isotopic to $\xi_n$ if 
$m \neq n$ by \fullref{nonisot}.
\end{proof}
%%%%%%%%%%%%%%%%%%%%%%%%%%%%%%%%%%%%%%%%%%%%%%%%%%%%%%%%%%%%%%%%%%%
\section{Ozsv\'ath--Szab\'o contact invariants}
In this section we give a brief overview of those properties of 
Heegaard Floer homology and of the related Ozsv\'ath--Szab\'o contact 
invariant which will be used in this article. 

Heegaard Floer homology is a family of functors
introduced by Ozsv\'ath and 
Szab\'o in \cite{O-Sz:2,O-Sz:1,O-Sz:3} which, in their simplest form, 
associate finitely generated Abelian groups 
$\widehat{HF}(Y, \mathfrak{t})$ to any closed connected\footnote{This
is the main deviation of the properties of
Heegaard Floer homology from the axioms of a topological quantum 
field theory. We thank Tom Mrowka for pointing out this issue.} 
oriented $\Spin^c$ $3$--manifold $(Y, \mathfrak{t})$, and  homomorphisms 
\[ \wwhat{F}_{W, \mathfrak{s}} \colon \widehat{HF}(Y_1, \mathfrak{t_1}) \to 
\widehat{HF}(Y_2, \mathfrak{t_2}) \]
to any oriented $\Spin^c$--cobordism $(W, \mathfrak{s})$ between 
two $\Spin^c$--manifolds $(Y_1, \mathfrak{t_1})$ and 
$(Y_2, \mathfrak{t_2})$ such that $\mathfrak{s}|_{Y_i}= \mathfrak{t_i}$.
If we do not need to specify the $\Spin^c$--structure on $W$ we write
\[ \wwhat{F}_W= \hspace{-10pt}\sum_{\scriptstyle {\mathfrak{s} \in 
\Spin^c(W)} \atop {\mathfrak{s}|_{Y_1}= \mathfrak{t}_1, \;
\mathfrak{s}|_{Y_2}= \mathfrak{t}_2}} 
\hspace{-10pt}\wwhat{F}_{W, \mathfrak{s}}. \] 
This notation makes sense because $\wwhat{F}_{W, \mathfrak{s}} \neq 0$ only
for finitely many $\Spin^c$--structures $\mathfrak{s}$.

A feature of Heegaard Floer homology is that, when $Y$ is a 
rational homology sphere, $\chi(\widehat{HF}(Y, \mathfrak{t}))=1$ for
 any $\mathfrak{t} \in \Spin^c(Y)$, where
the Euler characteristic is computed using a suitably defined 
$\Z / 2 \Z$--grading (see Ozsv\'ath--Szab\'o
\cite[Proposition~5.1]{O-Sz:2}). This implies 
that ${\rm rk} \widehat{HF}(Y, \mathfrak{t}) \geq 1$ for any 
$\Spin^c$--structure $\mathfrak{t}$.
\begin{dfn}
A rational homology sphere $Y$ is an {\em $L$--space} if
$ \widehat{HF}(Y, \mathfrak{t}) \cong \Z$ for any 
$\mathfrak{t} \in \Spin^c(Y)$.
\end{dfn}

A contact structure $\xi$ on a $3$--manifold $Y$ determines a 
$\Spin^c$--structure $\mathfrak{t}_{\xi}$ on $Y$ such that 
$c_1(\mathfrak{t}_{\xi})=c_1(\xi)$. To any contact manifold $(Y, \xi)$ we can 
associate an element $c(\xi) \in \widehat{HF}(-Y, \mathfrak{t}_{\xi})/ \pm 1$ 
which is an isotopy invariant of $\xi$, see \cite{O-Sz:cont}.
In the following we will always abuse the notation and consider 
$c(\xi)$ as an element of $\widehat{HF}(-Y, \mathfrak{t}_{\xi})$, 
although it is, strictly speaking, defined only up to sign. This
abuse does not lead to mistakes as long as we do not use the 
additive structure on $\widehat{HF}(-Y, \mathfrak{t}_{\xi})$.

\begin{thm}[Ozsv\'ath--Szab\'o {{\cite[Theorem~1.4 and
Theorem~1.5]{O-Sz:cont}}}]
If $(Y, \xi)$ is overtwisted, then $c(\xi)=0$. If $(Y, \xi)$ is Stein 
fillable, then $c(\xi)$ is a primitive element of 
$\widehat{HF}(-Y, \mathfrak{t}_{\xi})$.
\end{thm}
\begin{thm}[Ozsv\'ath--Szab\'o {{\cite[Theorem~4.2]{O-Sz:cont}}},
Lisca--Stipsicz {{\cite[Theorem 2.3]{lisca-stipsicz:3}}}]
\label{chirurgia}
Suppose that $(Y', \xi')$ is obtained from $(Y, \xi)$ by Legendrian 
surgery along a Legendrian link. Then we have
\[ \wwhat{F}_W(c(\xi'))= c(\xi) \] 
where $W$ is the cobordism induced by the surgery viewed as a 
cobordism from $-Y'$ to $-Y$.
\end{thm}  
We observe that it makes sense to write $\wwhat{F}_W(c(\xi'))$
because $\wwhat{F}_W$  descends to a well defined map
\[ \wwhat{F}_W \colon \widehat{HF}(-Y', \mathfrak{t}_{\xi'})/ \pm 1 \to 
\widehat{HF}(-Y, \mathfrak{t}_{\xi}) / \pm 1.
 \]
There is a also a more refined version of Heegaard Floer homology 
called { \em Heegaard Floer homology with twisted coefficients}; see 
Ozsv\'ath--Szab\'o \cite[Section~8]{O-Sz:2}.
Denote by $\Z[H^1(Y, \Z)]$ the group ring of the first cohomology 
group of $Y$ with integer coefficients. For any 
$\Z[H^1(Y, \Z)]$--module $M$ and for any compact connected
oriented $3$--manifold $Y$ endowed with a $\Spin^c$--structure 
$\mathfrak{t}$ one can define a Heegaard Floer homology 
group $\uline{\widehat{HF}}(Y, \mathfrak{t}; M)$.  
On this Heegaard Floer homology group there 
is a natural structure of module over $\Z[H^1(Y, \Z)]$ inherited 
by the coefficient group $M$. 

The contact invariant in the context of Heegaard Floer homology with 
twisted coefficients is denoted by $\underline{c}(\xi; M)$. It is well 
defined only up to multiplication by invertibles elements of 
$\Z[H^1(Y, \Z)]$, therefore it is, properly speaking, an element 
of the quotient 
$\uline{\widehat{HF}}(Y, \mathfrak{t}; \Z) / \Z[H^1(Y, \Z)]^*$.
If we consider $\Z[H^1(Y, \Z)]$ as a module over itself, we write
$\uline{\widehat{HF}}(Y, \mathfrak{t})$ for 
$\uline{\widehat{HF}}(Y, \mathfrak{t}; \Z[H^1(Y, \Z)])$, and
$\underline{c}(\xi)$ for $\underline{c}(\xi; \Z[H^1(Y, \Z)])$.
 The untwisted Heegaard Floer homology group 
$\widehat{HF}(Y, \mathfrak{t})$ can be seen in the theory with 
twisted coefficients as $\uline{\widehat{HF}}(Y, \mathfrak{t};
\Z)$, where $\Z$ is considered as a trivial 
$\Z[H^1(Y, \Z)]$--module, and $c(\xi)= \underline{c}(\xi; \Z)$.

\begin{rem}
If $Y$ is a rational homology sphere, then $\Z[H^1(Y, \Z)]= \Z$, 
therefore $\uline{\widehat{HF}}(Y, \mathfrak{t})=
\widehat{HF}(Y, \mathfrak{t})$ and $\underline{c}(\xi)= c(\xi)$.
\end{rem}

We say that a contact manifold $(Y, \, \xi)$ is 
{\em weakly symplectically fillable} if there is a symplectic manifold
$(X, \, \Omega)$ such that $Y= \partial X$ and $\Omega|_{\xi}>0$.
The Ozsv\'ath--Szab\'o contact invariant with twisted coefficients 
is non trivial for weakly fillable contact manifolds. 
More precisely, consider the $\Z[H^1(Y, \Z)]$--module $\Z[\R]$ 
generated over $\Z$ by the elements $T^r$, where $T$ is a formal 
variable and $r$ is any real number. A closed $2$--form $\omega$ on $Y$
 endows $\Z[\R]$ with the $H^1(Y, \Z)$--action $\gamma \cdot T = T^{\langle \gamma \cup [\omega], Y \rangle}$.
 Denote by $\uline{\widehat{HF}}(Y, \mathfrak{t}; [\omega])$ the 
Heegaard Floer homology group of $(Y, \mathfrak{t})$ with twisted 
coefficients in $\Z[\R]$ with the module structure defined by $\omega$.
\begin{thm}[Ozsv\'ath--Szab\'o {{\cite[Theorem~4.2]{O-Sz:genus}}}]
Let $(Y, \xi)$ be a weakly fillable
contact $3$--manifold, and let $(X, \Omega)$ be one of its weak 
fillings. Call $\omega= \Omega|_{Y}$, then $\underline{c}(\xi, [\omega])$ is a 
primitive element of $\uline{\widehat{HF}}(- Y, \mathfrak{t}_{\xi}; 
[\omega])$.
\end{thm}

\begin{rem} \label{wlog}
If $\omega$ is exact (in particular, if $Y$ is a rational homology 
sphere), then $\uline{\widehat{HF}}(Y, \mathfrak{t}_{\xi}; [\omega]) \cong
\widehat{HF}(Y, \mathfrak{t}_{\xi}) \otimes_{\Z} \Z[\R]$ and 
$\underline{c}(\xi, [\omega])= c(\xi) \otimes 1$.
\end{rem}

In the following we will use only the untwisted invariant because 
we are concerned only with rational homology spheres. We have 
made this excursion into Heegaard Floer homology with twisted 
coefficients only to show that we are not losing any generality
by considering only the untwisted invariant in this article.
%%%%%%%%%%%%%%%%%%%%%%%%%%%%%%%%%%%%%%%%%%%%%%%%%%%%%%%%%%%%%%%%%%
\section{Computation of the Ozsv\'ath--Szab\'o invariants}
We construct contact structures $\zeta_n$ on $M(r_1, r_2, r_3)$ for 
$n \geq 0$ so that we obtain $(M(r_1, r_2, r_3, r_4), \, \xi_n)$ from 
$(M(r_1, r_2, r_3), \, \zeta_n)$ by a sequence of negative Legendrian 
surgeries.

The contact manifold $(M \setminus V_4, \, \xi_n|_{M \setminus V_4})$ has infinite 
boundary slope for any $n \geq 0$. Let $\lambda$ be the unique tight contact 
structure with infinite boundary slope on the solid torus 
$D^2 \times S^1$. Then $(M(r_1, r_2, r_3), \, \zeta_n)$ is obtained 
 by gluing $(D^2 \times S^1, \lambda)$ to  $(M \setminus V_4, \, \xi_n|_{M \setminus V_4})$ by a map 
$\partial D^2 \times S^1 \to - \partial (M \setminus V_4)$ represented by the identity matrix in the 
bases described in \fullref{tre}.

$\partial D^2 \times S^1 \to - \partial (M \setminus V_4)$ 

\begin{prop} \label{homotopy}
All the contact structures $\xi_n$ are homotopic to $\xi_0$. All the 
contact structures $\zeta_n$ are homotopic to $\zeta_0$.
\end{prop}
\begin{proof}
All the contact structures $\xi_n$ coincide with 
$\wtilde{\xi}_+$ on $M_{c_1}$, and with 
$\wtilde{\xi}_-$ on $M_{c_2}$, therefore we need to 
show that they are homotopic relative to the boundary on 
$M \setminus (M_{c_1} \cup M_{c_2})$. We recall that 
$\xi_n|_{M \setminus (M_{c_1} \cup M_{c_2})}$ is a perturbation of 
$\alpha_{2n}(s_1,s_2)|_{T^2 \times [c_1,c_2]}$ near to the boundary in 
order to make it convex. If we choose the functions $\varphi_{2n}$ so that 
they coincide on $[-1, c_1+ \epsilon]$ and on $[c_2- \epsilon, 1]$, and do the 
perturbations in $[-1, c_1+ \frac{\epsilon}{2}] \cup [c_2- \frac{\epsilon}{2}, 1]$, 
we only need to show that all $\alpha_{2n}$ are homotopic through a 
homotopy which is constant in $[-1, c_1+ \epsilon] \cup [c_2- \epsilon, 1]$.

In order to construct such a homotopy, take a cut-off function $\beta$ 
which is $0$ on the union $[-1, c_1+ \frac{\epsilon}{2}] \cup [c_2- \frac{\epsilon}{2}, 1]$ and 
$1$ on $[c_1+ \epsilon, c_2- \epsilon]$. For any $n_0$ and $n_1$ we can take a number
 $K \gg 0$ which is big enough so that,
for any $\lambda \in [0,3]$, the kernel of the $1$--form $H_{\lambda}$ defined by
 \[\!  H_{\lambda}= \left \{ \begin{array}{l}  
\lambda K \beta(t)dt + \big ( \cos (\varphi_{2n_0}(t))dx+ \sin (\varphi_{2n_0}(t))dy \big ) 
\quad \text{if} \quad \lambda \in [0, 1 ]  \\ \\
K \beta(t)dt + (2- \lambda) \big ( \cos (\varphi_{2n_0}(t))dx+ \sin (\varphi_{2n_0}(t))dy 
\big ) + \\
(\lambda -1) \big ( \cos (\varphi_{2n_1}(t))dx+ \sin (\varphi_{2n_1}(t))dy \big ) \quad 
\text{if} \quad \lambda \in [1, 2] \\ \\
(3- \lambda) K \beta(t)dt + \big ( \cos (\varphi_{2n_1}(t))dx+ \sin (\varphi_{2n_1}(t))dy 
\big ) \! \! \quad \! \text{if} \! \! \quad \! \lambda \in [2, 3]
 \end{array} \right. \]
is a tangent $2$--plane field. A smoothing of $H_{\lambda}$ provides  
the wanted homotopy. The homotopy between $\zeta_{n_0}$ and $\zeta_{n_1}$ 
follows at once because $(N, \, \zeta_{n_0})$ and $(N, \, \zeta_{n_1})$ are 
obtained by modifying $(M, \, \xi_{n_0})$ and $(M, \, \xi_{n_1})$ in  
$M_{c_2}$ where the homotopy $H_{\lambda}$ is constant.
\end{proof}

\begin{prop} \label{surgery}
$(M(r_1, r_2, r_3, r_4), \, \xi_n)$ is obtained from $(M(r_1, r_2, r_3), \, \zeta_n)$ 
by a sequence of Legendrian surgeries.
\end{prop}
\begin{proof}
$(D^2 \times S^1, \lambda)$ is diffeomorphic to a standard neighbourhood of a 
Legendrian curve with twisting number $0$. Moreover the core of 
$D^2 \times S^1$ is isotopic in $M(r_1,r_2,r_3)$ to a regular fibre, 
therefore $(M \setminus V_4, \, \xi_n|_{M \setminus V_4})$ is the complement of the 
standard neighbourhood of a regular fibre with twisting number 
$0$ in $M(r_1,r_2,r_3)$. This implies that  $(M(r_1, r_2, r_3, r_4), \, \xi_n)$ is 
obtained from $(M(r_1, r_2, r_3), \, \zeta_n)$ by rational contact surgery.
 Since we have performed the surgery on a Legendrian curve with 
twisting number $0$, the contact surgery coefficient is equal to
the smooth surgery coefficient, which is $- \frac{1}{r_4} <0$. 
 By Ding--Geiges \cite[Proposition~3]{ding-geiges:2} any contact surgery with 
negative coefficient can be expanded into a sequence of Legendrian 
surgeries.
\end{proof}

\begin{prop} \label{zeta}
$(M(r_1, r_2, r_3), \, \zeta_0)$ is Stein fillable. $(M(r_1, r_2, r_3), \, \zeta_n)$
 is overtwisted for $n \geq 1$.
\end{prop}
\begin{figure}\centering\small
\begin{tabular}{ccc}
\includegraphics[width=5cm]{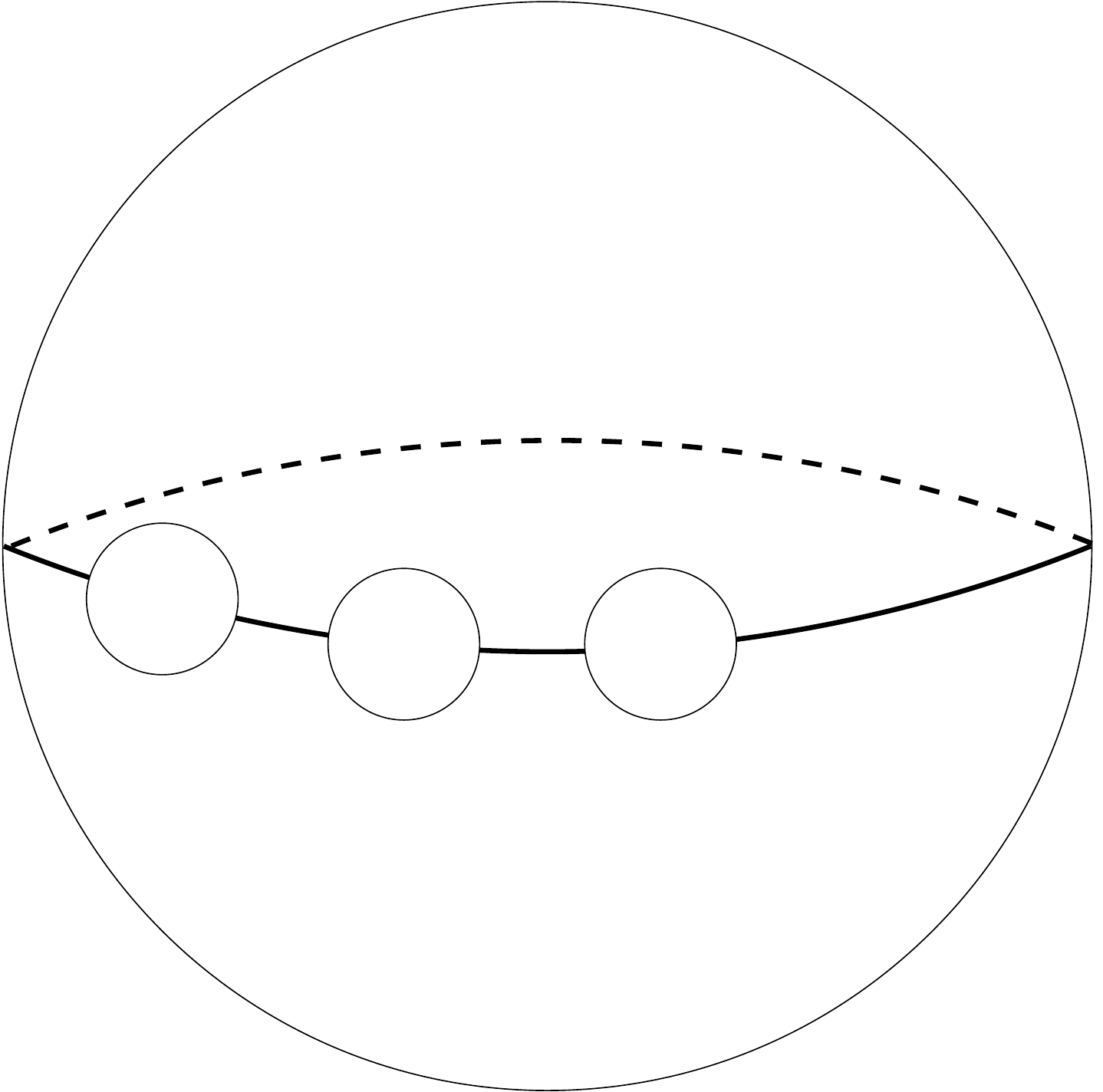}

& & 

\includegraphics[width=5cm]{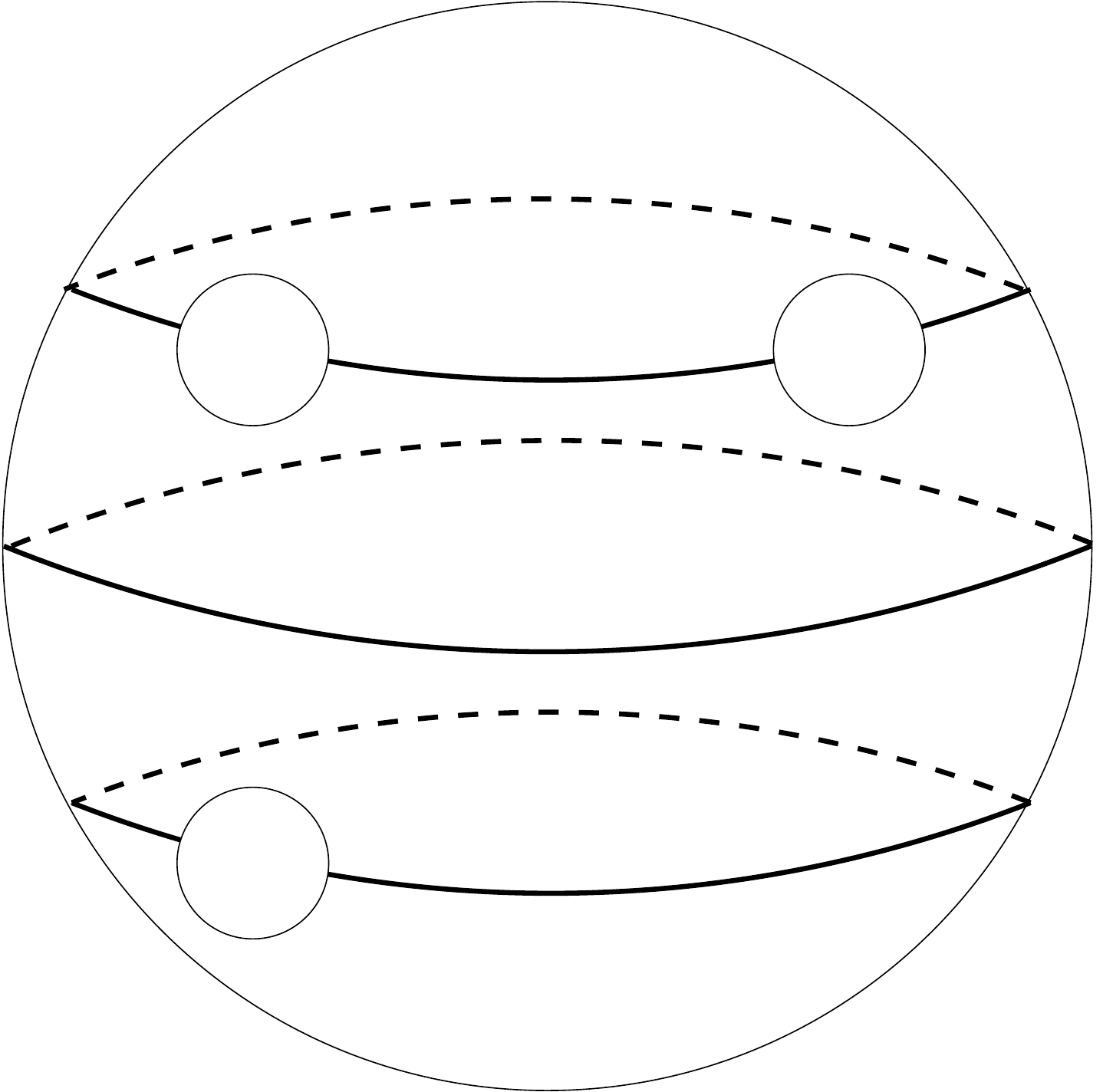} \\
(a) & & (b)

\end{tabular}
\caption{On the left the dividing set on a $\# \Gamma$--minimising 
section of the background of $(N, \, \zeta_0)$. On the right the 
dividing set on a $\# \Gamma$--minimising section of the background of 
$(N, \,\zeta_1)$.}
\label{background2.fig}
\end{figure}

\begin{proof}
The background of $(N, \, \zeta_n)$ is isomorphic to a $S^1$--invariant 
contact structure in $S' \times S^1$ where $S'$ is a three-punctured 
sphere. A convex $\# \Gamma$--minimising section $S_0'$ of the background
 of $(N, \, \zeta_n)$ is obtained by gluing a meridional disc of 
$(D^2 \times S^1, \lambda)$ to $S_0$ along the component of $\partial S_0$ corresponding 
to $V_4$. The dividing set of $S_0'$ consists of 
\begin{enumerate}
\item three arcs joining the boundary components of $\partial S_0'$ in 
pairs when $n=0$ (\fullref{background2.fig}(a)), or
\item two arcs joining two boundary components, one arc with both 
endpoints on the third boundary component, and $2n-1$ closed 
curves parallel to the third boundary component when $n \geq 1$ 
(\fullref{background2.fig}(b)). 
\end{enumerate}

We consider $\zeta_0$ first. It has been proved by Ghiggini, Lisca and
Stipsicz \cite{gls:1} 
that any tight contact structure on $N=M(r_1,r_2,r_3)$ which is
isomorphic to the background of $(N, \zeta_0)$ in the complement 
of tubular neighbourhoods of the singular fibres is Stein fillable 
independently of its restrictions to the neighbourhoods of the 
singular fibres, provided that the restrictions are tight. 

Now we consider $\zeta_n$ for $n \geq 1$.
The dividing arc with both endpoints on $- \partial (N \setminus V_3)$ produces a 
singular bypass on $S_0'$ by Honda \cite[Proposition 3.18]{honda:1}. 
By \cite[Lemma 3.15]{honda:1} attaching this bypass to 
$- \partial (N \setminus V_3)$ we thicken $V_3$ to $V_3'$ so that $- \partial (M \setminus V_3')$ 
has slope $0$. 

Take a point $p$ belonging to another dividing curve of $S_0'$, 
then $\{ p \} \times S^1$ is a Legendrian fibre with twisting number $0$ 
because $\zeta_n|_{N \setminus (V_1 \cup V_2 \cup V_3)}$ is $S^1$--invariant by 
Honda \cite[Section 4.3]{honda:2}.
Applying the Imbalance principle \cite[Proposition 3.17]{honda:1}
we use this curve to  find a vertical bypass attached to 
$\partial (N \setminus V_3')$. The attachment of this bypass gives a further 
thickening of $V_3'$ to $V_3''$ so that $- \partial (N \setminus V_3'')$ has 
infinite boundary slope again. By 
\cite[Proposition~4.16]{honda:1}  there is a standard torus with 
slope $-r_3$ in $V_3'' \setminus V_3$. This torus produces an overtwisted disc 
in $(N, \zeta_n)$. 
\end{proof}

\begin{cor} \label{xi}
$(M(r_1, r_2, r_3, r_4), \, \xi_0)$ is Stein fillable.
\end{cor}

\begin{thm} \label{annullamento}
$c(\xi_n)=0$ for $n \geq 1$.
\end{thm}
\begin{proof}
 The contact structures $\xi_n$ 
are homotopic for $n \geq 0$ by \fullref{homotopy}, in particular they 
determine the same $\Spin^c$--structure on $M$.
The same is true for the contact structures $\zeta_n$ on $N$.
Let $\mathfrak{t}_{\xi}$ denote the $\Spin^c$--structure on $M$ 
determined by the $\xi_n$'s, and let $\mathfrak{t}_{\zeta}$ denote the 
$\Spin^c$--structure on $N$ determined by the $\zeta_n$'s. 

The surgery links from \fullref{surgery} are all smoothly isotopic
independently of the contact structure, therefore they determine the
same smooth cobordism $W$ from $-M$ to $-N$. This implies that for any
$n \geq 0$ we have $\wwhat{F}_W(c(\xi_n))= c(\zeta_n)$. Both $- M$ and
$- N$ are $L$--spaces by Ozsv\'ath--Szab\'o
\cite[Lemma~2.6]{O-Sz:plumbing}, hence $\widehat{HF}(-M,
\mathfrak{t}_{\xi})$ is generated by $c(\xi_0)$ and $\widehat{HF}(-N,
\mathfrak{t}_{\zeta})$ is generated by $c(\zeta_0)$ because $\xi_0$
and $\zeta_0$ are Stein fillable by \fullref{xi} and
\fullref{zeta}. This implies that $\wwhat{F}_W$ is injective. 
$c(\zeta_n)=0$ when $n \geq 1$ because $\zeta_n$ is overtwisted by
\fullref{zeta}, therefore the injectivity of $\wwhat{F}_W$ implies
that $c(\xi_n)=0$.
\end{proof}

$c(\zeta_n)=0$ when 
$n \geq 1$ because $\zeta_n$ is overtwisted by \ref{zeta}, therefore

\begin{cor}\label{nonfill}
$(M(r_1, r_2, r_3, r_4), \, \xi_n)$ is not weakly symplectically fillable
 for $n \geq 1$.
\end{cor}
%%%%%%%%%%%%%%%%%%%%%%%%%%%%%%%%%%%%%%%%%%%%%%%%%%%%%%%%%%%%%%%%%
\section{Further examples}
More examples of tight contact manifolds with trivial contact
invariants can be constructed by Legendrian surgery on 
$(M(r_1,r_2,r_3), \, \xi_n))$ with the help of the following tightness 
criterion, which was implicitly proved in Hofer's work on the Weinstein 
conjecture for overtwisted contact structures \cite{hofer:1}.
\begin{prop} \label{hofer}
Let $(Y, \, \xi)$ be a hypertight contact manifold. Then any contact
manifold $(Y', \, \xi')$ obtained by Legendrian surgery on 
$(Y, \, \xi)$ is tight.
\end{prop}
\begin{proof}[Sketch of the proof] Pick a contact form $\alpha$ on $Y$ 
whose Reeb flow has no contractible periodic orbits. The Legendrian 
surgery defines an exact symplectic cobordism $(W, d \lambda)$ from 
$(Y, \, \xi)$ to $(Y', \, \xi')$ such that $\lambda|_Y= \alpha$. Extend $W$ to a 
non compact exact symplectic manifold 
$(\wwhat{W}, \, d \what{\lambda})$ by gluing 
$(Y \times (- \infty ,0], \, d(e^t \alpha))$ to $Y \subset W$.

Suppose that $(Y', \xi')$ is overtwisted, then there is a $2$--sphere
$S$ embedded in $Y'$ whose characteristic foliation contains a 
closed periodic orbit and has an elliptic point in each connected 
component of the complement of the periodic orbit as unique 
singularities. Now we proceed as in \cite{hofer:1}: we start 
filling $S$ by a Bishop family of holomorphic discs originating 
from an elliptic singularity of the characteristic foliation of
 $S$. Since $S$ cannot be filled completely because of the periodic 
orbit in the characteristic foliation,  we end up with a 
contractible periodic orbit of the Reeb flow associated to $\alpha$ 
exactly like in the case of the symplectisation, contradicting our 
hypothesis.
\end{proof}

Isotope the singular fibres $F_i$ to Legendrian curves with twisting 
number $-1$ in $(M(r_1, r_2, r_3, r_4), \, \xi_n)$. Denote by 
$(M(r_1'', r_2'', r_3'', r_4''), \, \xi_n'')$ a contact manifold obtained 
by contact surgery on the singular fibres $F_i$ with smooth
surgery coefficient $r_i'$. From well known properties of rational
surgery and Seifert invariants we have
$\smash{- \frac{1}{r_i''}= -\frac{1}{r_1} - \frac{1}{r_1'}}$,
therefore $\smash{r_i''= \frac{r_ir_i'}{r_i+r_i'}}$.
 Observe that, in general, $\xi_n''$ is not uniquely determined by the 
surgery coefficients, however we will not care about this fact 
because our argument is independent of the choices determining 
$\xi_n''$.

If $r_i' < -1$ by Ding--Geiges \cite[Proposition~3]{ding-geiges:2} any contact 
surgery on $F_i$ with smooth surgery coefficient $r_i'$ can be 
realised as a Legendrian surgery on a Legendrian link in 
$(M(r_1, r_2, r_3, r_4), \, \xi_n)$, and $\xi_n''$ is determined by the rotation
numbers of the components of the Legendrian surgery link. By 
\fullref{hofer} the resulting contact manifold 
$(M(r_1'', r_2'', r_3'', r_4''), \, \xi_n'')$ is tight. 

The proof of \fullref{annullamento} goes through for this more general family of 
contact manifolds, therefore we have the following theorem.
\begin{thm}
If $r_i' < -1$ and $n>0$, for any choice of the rotation numbers of 
the components of the Legendrian surgery link defining $\xi_n''$ we 
have $c(\xi_n'')=0$. 
\end{thm}
Most of the contact structures $\xi_n''$ are virtually overtwisted. 
To prove this we apply \fullref{thickening} to 
$(M(r_1'', r_2'', r_3'', r_4''), \, \xi_n'')$ in order to find tubular 
neighbourhoods $V_i$ of the singular fibres, and observe that the 
restriction of $\xi_n''$
to $V_i$ may contain basic slices with different signs, depending on
the choices made in the contact surgery. In those cases
$\xi_n''|_{V_i}$, and therefore $\xi_n''$, is virtually overtwisted.

An alternative way to prove that most contact structures $\xi_n''$ are
virtually overtwisted is by applying Gompf's trick 
\cite[Proposition~5.1]{gompf:1} to the Legendrian surgery diagrams 
for $(M(r_1'', r_2'', r_3'', r_4''), \, \xi_n'')$ having components with
non maximal rotation number in order to find an overtwisted disc in 
a finite cover. 

The proof of \fullref{torsione} does not extend immediately 
to $\xi_n''$ because the contact structures considered by 
Honda--Kazez--Mati\'c \cite{honda-kazez-matic:2} are  universally tight, therefore we are 
not able to prove that they are all distinct. However we 
believe that $\xi_n''$ and $\xi_m''$ are non isomorphic if $m \neq n$, and 
that two contact structures on $M(r_1'',r_2'',r_3'',r_4'')$ obtained by 
Legendrian surgery on $(M(r_1,r_2,r_3,r_4), \, \xi_n)$ are non isotopic if 
they are obtained by surgery on Legendrian links with different 
rotation numbers.

\section{Final considerations}
Giroux has introduced the a topological invariant for 
$3$--dimensional contact manifolds defined as follows.
\begin{dfn}  
Let $\xi$ be a contact structure on a $3$--manifold $Y$. The 
{\em torsion} of $(Y, \, \xi)$ is the supremum of the integers $n \geq 1$
for which there exists a contact embedding
\[ 
\big ( T^2 \times [0,1],  \, \ker (\cos (2n \pi z)dx - \sin (2n \pi z)dy) 
\big ) \hookrightarrow (Y, \, \xi).
\] 
 We declare the torsion of 
$(Y, \, \xi)$ to be $0$ if no such an embedding exists.
\end{dfn}
We denote the torsion of $(Y, \, \xi)$ by ${\rm Tor}(Y, \, \xi)$.
One can deduce from \fullref{nonisom} that 
${\rm Tor}(M, \, \xi_n)=n$, therefore \fullref{nonfill} adds further evidence to the
following conjecture, which the author learnt from Eliashberg.
\begin{conj}\label{yashascon}
If ${\rm Tor}(Y, \, \xi) >0$ then $(Y, \, \xi)$ is not strongly 
symplectically fillable.
\end{conj}
It seems also plausible that the following stronger statement holds.
\begin{conj} \label{miacon}
If ${\rm Tor}(Y, \, \xi) >0$ then $c(\xi)=0$.
\end{conj}
\fullref{miacon} implies \fullref{yashascon} because 
a strongly symplectically fillable contact manifold has non trivial 
contact invariant by \fullref{wlog}.
On the other hand there are 
contact manifolds with positive torsion which are weakly fillable, 
therefore they have non trivial 
twisted contact invariant.
 
We are able to define another family of universally tight contact structures
$\eta_n$ on $M(r_1, r_2, r_3, r_4)$ for $n \geq 0$ coinciding with
$\xi_+$ on $M'(r_1,r_2)$ and on $M'(r_3,r_4)$, and  with $\alpha_{2n+1}(s_1,s_2)$ 
on $T^2 \times [-1,1]$.
We observe that ${\rm Tor}(M, \eta_n)=n$. Since all the example of
tight contact structures with trivial Ozsv\'ath--Szab\'o invariants 
we know at present have positive torsion, it would be interesting 
to compute $c(\eta_0)$. Unfortunately the strategy adopted in this 
article fails for computing $c(\eta_n)$, because we cannot show that
the $\eta_n$'s are homotopic to Stein fillable contact structures.
%%%%%%%%%%%%%%%%%%%%%%%%%%%%%%%%%%%%%%%%%%%%%%%%%%%%%%%%%%%%%%%%%%

\bibliographystyle{gtart}
\bibliography{link}

\end{document}